\documentclass{amsart}

\input xy
\xyoption{all}

\newtheorem{theorem}{Theorem}[section]
\newtheorem{proposition}[theorem]{Proposition}
\newtheorem{lemma}[theorem]{Lemma}
\newtheorem{corollary}[theorem]{Corollary}
\newtheorem{fact}[theorem]{Fact}

\theoremstyle{definition}
\newtheorem{definition}[theorem]{Definition}

\theoremstyle{remark}
\newtheorem{remark}[theorem]{Remark}

\numberwithin{equation}{section}

\def\acl{\operatorname{acl}}

\def\dim{\operatorname{dim}}
\def\tp{\operatorname{tp}}

\def\id{\operatorname{id}}

\def\loc{\operatorname{loc}}

\begin{document}

\title{An essentially saturated surface not of K\"ahler-type}

\author{Rahim Moosa}
\address{University of Waterloo \\
Department of Pure Mathematics MC 5049\\
200 University Avenue West\\
Waterloo, Ontario N2L 3G1 \\
Canada}
\thanks{Rahim Moosa was partially supported by NSERC}
\email{rmoosa@math.uwaterloo.ca}

\author{Ruxandra Moraru}
\address{University of Waterloo \\
Department of Pure Mathematics MC 5170\\
200 University Avenue West\\
Waterloo, Ontario N2L 3G1 \\
Canada}
\thanks{Ruxandra Moraru was partially supported by NSERC}
\email{moraru@math.uwaterloo.ca}

\author{Matei Toma}
\address{IECN,
Nancy-Universit\'e,  CNRS, INRIA\\
Boulevard des 
Aiguillettes B.P. 239, F-54506 Vandoeuvre-l\`es-Nancy Cedex\\
France \\and Mathematical Institute of the Romanian Academy}
\email{toma@iecn.u-nancy.fr}

\date{December 19, 2007}

\begin{abstract}
It is shown that if $X$ is an Inoue surface of type $S_M$ then the irreducible components of the Douady space of $X^n$ are compact, for all $n\geq 0$.
This gives an example, asked for in~\cite{sat}, of an {\em essentially saturated} compact complex manifold (in the sense of model theory) that is not of K\"ahler-type.
Among the known compact complex surfaces without curves, these are the only examples.
\end{abstract}

\maketitle

\section{Introduction}
\label{introduction}
\noindent
A (reduced) compact complex-analytic space $X$ can be viewed as a first-order structure, in the sense of mathematical logic, by equipping $X$ with a predicate 
symbol $P_A$ for each complex-analytic subset $A\subseteq X^n$, for all $n\geq 0$.
We denote this structure by $\mathcal{A}(X)$.
Model-theory is interested in the {\em definable sets} of $\mathcal{A}(X)$: the subsets of $X^n$, for various $n$, obtained from the complex-analytic sets by taking intersections, complements, fibres of co-ordinate projections, and images of co-ordinate projections.
Zilber~\cite{zilber93} showed that this structure is ``tame'' in that it admits {\em quantifier elimination} (every definable set is a finite boolean combination of complex-analytic sets) and a certain model-theoretic rank ({\em Morley rank}) is finite valued on definable sets.
Motivated by model-theoretic considerations, the first author introduced in~\cite{sat} the notion of an {\em essentially saturated} compact complex-analytic space, namely, one for which there exists a countable subcollection of the predicates $P_A$ from which one can define all the definable sets of $\mathcal{A}(X)$.
The main result in~\cite{sat}, slightly reformulated, is the following geometric characterisation:
{\em A compact complex-analytic space $X$ is essentially saturated if and only if, for all $n\geq 0$, every irreducible complex-analytic subset of $X^n$ lives in an irreducible component of the Douady space of $X^n$ that is compact.}
Recall that the Douady space is the analytic analogue of the Hilbert scheme; it parameterises all compact complex-analytic subspaces of $X$ (see Section~\ref{preliminaries} below for some details).

Every holomorphic image of a compact K\"ahler manifold is essentially saturated (these are the {\em K\"ahler-type} spaces introduced by Fujiki in~\cite{fujiki78}).
The first author asked in~\cite{sat} for an example of an essentially saturated space that is not of K\"ahler-type.
Akira Fujiki, in private communication, suggested that we consider the surfaces of type $S_M$ constructed by Inoue in~\cite{inoue74}.
The purpose of this note is to show that these surfaces are indeed examples of essentially saturated spaces that are not of K\"ahler-type.
In fact, among the known compact complex surfaces of non K\"ahler-type without curves, these are the only examples.
A key element in our proof is a model-theoretic classification, due to Pillay and Scanlon~\cite{pillayscanlon2000}, of compact complex manifolds having no proper infinite complex-analytic subsets (see Fact~\ref{trichotomy} below).

\smallskip

For the rest of this paper, by a {\em complex variety} we will mean a reduced and irreducible compact complex-analytic space.
A {\em subvariety} will mean an irreducible complex-analytic subset.

\bigskip
\section{Preliminaries}
\label{preliminaries}

\noindent
Let $X$ be any complex-analytic space.
There exist a complex analytic space $D(X)$, called the {\em Douady space} of $X$, and a complex-analytic subspace $Z(X)\subseteq D(X)\times X$, called the {\em universal family} of $X$, such that
\begin{itemize}
\item[(a)]
the projection $Z(X)\to D(X)$ is a flat and proper surjection, and
\item[(b)]
(Universal Property)
if $S$ is any complex-analytic space and $G$ is any complex-analytic subspace of $S\times X$ with $G\to S$ a flat and proper surjection, then there exists a unique holomorphic map (the {\em Douady map}) $\phi:S\to D(X)$ inducing a canonical isomorphism $G\simeq S\times _{D(X)} Z(X)$.
\end{itemize}
In particular, given a complex-analytic subset $A\subseteq X$ there is a unique point $d\in D(X)$ such that $A$ is the fibre, $Z(X)_d$, of $Z(X) \rightarrow D(X)$ at $d$.
This point, often denoted by $[A]$, is called the {\em Douady point} of $A$ in $X$.
The Douady space was constructed by Douady in~\cite{douady66} and shown to have countably many irreducible components by Fujiki in~\cite{fujiki79}.
A more detailed discussion of Douady spaces can be found in~\cite{campanapeternell94}.

It is not necessarily the case that if $X$ is a compact complex variety then the irreducible components of $D(X)$ are again compact.
Indeed, as explained in the introduction, the compactness of the components of the Douady space turns out to be model-theoretically very significant:
Essential saturation is equivalent to asking that every subvariety of $X^n$ live in an irreducible component of the Douady space of $X^n$ that is compact.
Here, given an irreducible component $D$ of $D(X^n)$, and a subvariety $A\subseteq X^n$, by ``$A$ lives in $D$'' we mean that $[A]\in D$.

Fujiki showed in~\cite{fujiki78} that the irreducible components of the Douady spaces of K\"ahler-type spaces (these are the holomorphic images of compact K\"ahler manifolds) are compact.
Since being of K\"ahler-type is preserved under taking cartesian products, this implies that K\"ahler-type spaces are essentially saturated.
A Hopf surface $H$ is an example of a compact complex manifold that is not essentially saturated.
While the components of $D(H)$ itself are compact,
$D(H\times H)$ has non-compact components coming from families of graphs of automorphisms of $H$ (see~\cite{campanapeternell94} for details).

We will show that the Inoue surfaces of type $S_M$ introduced by Inoue in~\cite{inoue74} are essentially saturated but not of K\"ahler-type.
Instead of recalling the construction of these surfaces, we collect together in the following fact those properties of these surfaces that will be relevant to our argument.

\begin{fact}
\label{inouefacts}
Suppose $X$ is an Inoue surface of type $S_M$. Then
\begin{itemize}
\item[(a)]
$X$ is a smooth compact complex surface containing no curves.
\item[(b)]
$H^i(X,T_X)=0$ for $i=0,1$.
\item[(c)]
$X$ is not of K\"ahler-type.
\item[(d)]
Any unramified covering of $X$ satisfies properties (a)-(c).
\end{itemize}
\end{fact}

\begin{proof}
Parts~(a) and~(b) are Proposition~2 of Inoue's original paper~\cite{inoue74}.
Part~(c) follows from the fact that the odd Betti numbers of every compact complex variety of K\"ahler-type are even~(\cite{fujiki83a}), while Inoue surfaces have first Betti number equal to one~(\cite{inoue74}).

Part~(d) follows from the fact that an unramified covering of an Inoue surface of type $S_M$ is again an Inoue surface of type $S_M$.
To see this, recall first of all that the Inoue surfaces are known to be exactly those smooth compact complex surfaces that have no curves, have first Betti number equal to one, and have second Betti number equal to zero
(cf. \cite{inoue74, LYZ94, Tel94}).
All of these properties are preserved under taking unramified coverings (this is Lemma~1 of~\cite{inoue74}).
Moreover, it follows from Inoue's constructions (see also p. 586 of ~\cite{brunella97}), that among the Inoue surfaces, those of type $S_M$ are distinguished as the only ones admitting two {\em holomorphic foliations}, i.e., whose tangent bundles admit two holomorphic subbundles.
As this property is also preserved under taking unramified coverings, we see that an unramified covering of an Inoue surface of type $S_M$ is again an Inoue surface of type $S_M$.
\end{proof}

We will eventually establish the essential saturation of Inoue surfaces of type $S_M$ by proving that the subvarieties of $X^n$ must be of a very special form: they are either ``degenerate'' (defined below) or live in a zero-dimensional component of the Douady space of $X^n$.
Lemma~\ref{isolated-degenerate-es} below shows that that will suffice.

\begin{definition}
Suppose $X$ is a compact complex variety.
A subvariety $A\subseteq X^n$ is {\em degenerate} if for some co-ordinate projection $p:X^n\to X$, $p(A)$ is a point.
\end{definition}

\begin{lemma}
\label{degeneratecomponents}
Degeneracy is preserved in components of the Douady space.
More precisely, suppose $A\subseteq X^n$ is a degenerate subvariety that lives in the irreducible component $D$ of $D(X^n)$ and let $Z$ be the restriction of $Z(X^n)$ to $D$.
Then $Z_d$ is degenerate for all $d\in D$.
\end{lemma}

\begin{proof}
Denote by $p:X^n\to X$ the first projection and by $\pi: D\times X^n\to D\times X$ the induced projection on $D\times X^n$. 
Up to a permutation of the co-ordinates we may suppose 
that $p(A)=a\in X$ is a point.
We will show in fact that $p(Z_d)$ is a point for all $d\in D$.
Note that the only obstruction is that $\pi(Z)\to D$ need not be flat.

First of all note that as $Z$ has a reduced and irreducible fibre, by flatness its general fibres are reduced and irreducible, as is $Z$ itself.
So $\pi(Z)$ and its general fibres over $D$ are irreducible.
But $\pi(Z)_{[A]}=\{a\}$.
Hence the general fibres of $\pi(Z)$ over $D$ must be points.
That is, $p(Z_d)$ is a point for general $d\in D$.

Let $D'$ be the set of points $d\in D$ such that $\dim p(Z_d)>0$.
We have shown that $D'$ is contained in a proper complex-analytic subset.
Now we claim that $D'$ is in fact empty.
Toward a contradiction, let $d'\in D'$.
Let $U$ be a sufficiently small open neighbourhood of $d'$ in $D$ and let $C$ be a complex-analytic curve in $U$ which passes through $d'$ and such that for general $c\in C$, $Z_c$ is irreducible and $p(Z_c)$ is a point
(this is possible as for general $d\in D$, $Z_d$ is irreducible and $p(Z_d)$ is a point).
Let $Z_C$ be the restriction of $Z$ to $C$.
By flatness, $Z_C$, and hence $\pi(Z_C)$, is irreducible.
Moreover the general fibres of $\pi(Z_C)\to C$ are points and so $\dim \pi(Z_C)=1$.
This contradicts the fact that the fibre of $\pi(Z_C)$ over $d'$ has positive dimension.
\end{proof}

\begin{lemma}
\label{isolated-degenerate-es}
Suppose $X$ is a compact complex variety with the property that every subvariety of $X^n$, for all $n\geq 0$, is either degenerate or lives in a zero-dimensional component of the Douady space of $X^n$.
Then $X$ is essentially saturated.
\end{lemma}

\begin{proof}
First of all, as explained in the introduction, it suffices to show that every subvariety of $X^n$ lives in an irreducible component of $D(X^n)$ that is compact.
We do this by induction on $n\geq 0$.
Cleary, the zero-dimensional components of $D(X^n)$ are compact, and so we focus on the degenerate subvarieties.
For $n=0$ there is nothing to prove.
For $n=1$ note that the degenerate subvarieties are just the points of $X$, and that $X$ itself is the (compact) component of $D(X)$ parametrising these points.

Now suppose $A$ is a degenerate subvariety of $X^n$, for some $n>1$, $D$ is the irreducible component of $D(X^n)$ in which $A$ lives, and $Z\subset D\times X^n$ is the universal family restricted to $D$.
As discussed in Lemma~\ref{degeneratecomponents}, after possibly permuting co-ordinates, for every $d\in D$, $Z_d$ is of the form $\{a\}\times B$ for some subvariety $B$ of $X^{n-1}$.
For each irreducible component $E$ of $D(X^{n-1})$ containing a subvariety of $X^{n-1}$, let $E'=X\times E$ and let $Z'\subset E'\times X^n$ be the complex-analytic subspace with $Z'_{(a,e)}=\{a\}\times Z(X^{n-1})_e$ for all $a\in X$ and $e\in E$.
Note that $E'$ is compact by induction, and $Z'\to E'$ is flat.
We have corresponding (injective) Douady maps $\phi_E:E'\to D(X^n)$ such that $Z(X^n)_{\phi_E(a,e)}=\{a\}\times Z(X^{n-1})_e$.
The $\phi_E(E')$'s must cover $D$ since every fibre above $D$ is of this form.
As there are only countably many  $\phi_E(E')$'s and each one is an irreducible complex-analytic subset of $D(X^n)$, it follows that $D$ must be equal to some $\phi_E(E')$, and thus be compact.
\end{proof}

\bigskip

\section{Trivial strongly minimal compact complex varieties}

\noindent
In this section we point out that if $X$ is a ``trivial strongly minimal" compact complex variety (explained below) with the property that every subvariety of $X\times X$ is either degenerate or lives in a zero-dimensional component of $D(X\times X)$, then the same is true of $X^n$ for all $n>2$.
In particular, such varieties will be essentially saturated.
This is not a very surprising result, since, as explained below, triviality says that all relations are essentially binary.

While we will use model-theoretic language freely in this section, we will try to give geometric formulations of the ideas involved and the results obtained.
We suggest~\cite{markerbook} as a general reference for model theory, and~\cite{moosa-ccs} for the model theory of compact complex varieties.

Let us begin by describing what the abstract notions of ``strong minimality'' and ``triviality'' amount to for compact complex varieties.
{\em Strong minimality} just means that $X$ has no proper infinite complex-analytic subsets.
Hence, for example, any irreducible compact complex surface without curves is strongly minimal.
{\em Triviality} is the following condition.
Suppose $n> 0$ and $A\subseteq X^{n+1}$ is a complex-analytic subset such that projection onto the first $n$ co-ordinates, $X^{n+1}\to X^n$, restricts to a generically finite-to-one map on $A$.
Then there must exist some $i\leq n$ such that if $A_i\subseteq X^2$ denotes the image of $A$ under the co-ordinate projection $(x_1,\dots,x_{n+1})\mapsto(x_i,x_{n+1})$, then the projection onto the first co-ordinate, $X^2\to X$, restricts to a generically finite-to-one map on $A_i$.

The following characterisation of non-trivial strongly minimal compact complex manifolds is a manifestation of the ``Zilber Trichotomy'' in this context.
It is due originally to Scanlon~\cite{scanlon2000} and appears as Proposition~5.1 of~\cite{pillayscanlon2000}.

\begin{fact}
\label{trichotomy}
If $X$ is a non-trivial strongly minimal compact complex manifold, then $X$ is either a complex torus or a projective curve.
In particular, every strongly minimal compact complex manifold that is not of K\"ahler-type is trivial.
\label{pillayscanlonfact}
\end{fact}

For example, since the Inoue surfaces of type $S_M$ have no divisors and are not of K\"ahler-type (see Fact~\ref{inouefacts} above), they are trivial strongly minimal.

The following fact about trivial strongly minimal sets in general is an easy and well-known model-theoretic consequence of the definitions.
If the reader unfamiliar with model theory is willing to accept Corollary~\ref{geom-trivial=binary} below, then he or she can skip this lemma and go directly to Proposition~\ref{reductionto2}.

\begin{lemma}
\label{trivial=binary}
Suppose $X$ is a trivial strongly minimal set in some saturated model of a complete stable theory.
Given $a=(a_1,\dots,a_n)\in X^n$, the partial $n$-type $\Sigma(x_1,\dots,x_n):=\displaystyle \bigcup_{i,j\leq n}\tp(a_i,a_j)$ has only finitely many completions.
That is, there exist finitely many complete $n$-types $q_1(x),\dots,q_N(x)$ such that given $b\in X^n$, if $\displaystyle b\models\Sigma$ then $b\models q_\ell$ for some $\ell\leq N$.
\end{lemma}

\begin{proof}
Without loss of generality we may assume that for some $r\leq n$, $\{a_1,\dots,a_r\}$ is an $\acl$-basis for $\{a_1,\dots,a_n\}$.
Suppose $b=(b_1,\dots,b_n)\models\Sigma$.
By triviality, any $\acl$-dependence among $\{b_1,\dots,b_r\}$ would be witnessed by a pair of elements in that set.
Since $\Sigma$ forces all pairs from the first $r$ co-ordinates to be $\acl$-independent, it follows that $\{b_1,\dots,b_r\}$ must be an $\acl$-independent set.
If $r=n$ then in fact $b$ and $a$ are generic tuples and so $b\models\tp(a)$.
That is, $\Sigma$ is complete and we are done.
Hence we may sssume that $r<n$.

Now for each $i=1,\dots,n-r$, $a_{r+i}\in\acl(a_1,\dots,a_r)$.
By triviality there is a $j_i\leq r$ such that $a_{r+i}\in\acl(a_{j_i})$.
Let $\phi_i(x_{j_i},x_{r+i})$ be a formula witnessing this.
So the set defined by $\phi_i(a_{j_i},x_{r+i})$ is finite and contains $a_{r+i}$.
Let $q_1,\dots,q_N$ be the set of all complete $n$-types of the form $\tp(a_1,\dots,a_r,c_1,\dots,c_{n-r})$ where each $c_i$ is in the set defined by $\phi_i(a_{j_i},x_{r+i})$.
Since $b$ realises $\Sigma$, $\models\phi_i(b_{j_i},b_{r+i})$ for all $i=1,\dots,n-r$.
Hence, if we let $f$ be an automorphism taking $(b_1,\dots,b_r)$ to $(a_1,\dots,a_r)$, then for each $i=1,\dots,n-r$ we have that $f(b_{r+i})$ is in the set defined by $\phi_i(a_{j_i},x_{r+i})$.
So $\tp(b)=\tp\big(a_1,\dots,a_r,f(b_{r+1}),\dots,f(b_n)\big)=q_\ell$ for some $\ell\leq N$.
\end{proof}

The following corollary gives the geometric content of Lemma~\ref{trivial=binary} specialised to compact complex varieties.

\begin{corollary}
\label{geom-trivial=binary}
Suppose $X$ is a trivial strongly minimal compact complex variety.
Given a subvariety $A\subseteq X^n$ there exist only finitely many other subvarieties having the same projections to $X\times X$.
More precisely, there exist
subvarieties $B_1,\dots,B_N\subseteq X^n$ such that if $B\subseteq X^n$ is a subvariety for which $\pi(A)=\pi(B)$ for all co-ordinate projections $\pi:X^n\to X^2$, then $B=B_\ell$ for some $\ell\leq N$.
\end{corollary}

\begin{proof}
We work in a saturated elementary extension $\mathcal{A}(X)'$ of $\mathcal{A}(X)$, which acts as ``universal domain'' (in the sense of Weil) for the geometry of the complex-analytic subsets of $X$ and its cartesian powers.
A key point is that the complete $n$-types in $\mathcal{A}(X)'$ (over the emptyset) are in one-to-one correspondence with subvarieties of $X^n$; every complete $n$-type is the generic type of a subvariety of $X^n$ called its {\em locus}.

Let $a$ realise the generic type of $A$ in $\mathcal{A}(X)'$.
Apply Lemma~\ref{trivial=binary} to $a$ to obtain complete types $q_1,\dots,q_N$.
For each $\ell\leq N$, let $B_\ell=\loc(q_\ell)$ be the locus of $q_\ell$.
Now suppose $B$ is as in the statement of the corollary and let $b$ realise the generic type of $B$.
Note that for each co-ordinate projection $\pi:X^n\to X^2$, $\pi(a)$ realises the generic type of $\pi(A)$ and $\pi(b)$ realises the generic type of $\pi(B)$.
So the assumption on $B$ says that $\displaystyle b\models\bigcup_{i,j\leq n}\tp(a_i,a_j)$.
Hence, by the conclusion of Lemma~\ref{trivial=binary}, $b\models q_\ell$ for some $\ell\leq N$.
It follows that $B=B_\ell$, as desired.
\end{proof}

Here is the main conclusion of this section.

\begin{proposition}
\label{reductionto2}
Suppose $X$ is a trivial strongly minimal compact complex variety with the property that every subvariety of $X\times X$ is either degenerate or lives in a zero-dimensional component of $D(X\times X)$.
Then the same holds of $X^n$, for each $n>2$.
In particular, $X$ is essentially saturated.
\end{proposition}

\begin{proof}
Note that every subvariety of $X$ is either degenerate or lives in a zero-dimensional component of $D(X)$: this follows from strong minimality as the only subvarietries of $X$ are points or $X$ itself.

Let us fix $n>2$ and a subvariety $Y\subseteq X^n$.
Assume $Y$ is not degenerate and does not live in a zero-dimensional component of $D(X^n)$, and seek a contradiction.
Let $D$ be the component of the Douady space of $X^n$ in which $Y$ lives and let $Z\subseteq D\times X^n$ be the restriction of the universal family to $D$.
Let $E$ be a proper complex-analytic subset of $D$ such that for all $d\in D\setminus E$, $Z_d$ is reduced and irreducible.
Since $Y$ is not degenerate, none of these $Z_d$'s are degenerate (cf. Lemma~\ref{degeneratecomponents}).
Hence $\pi(Z_d)$ is non-degenerate for each co-ordinate projection $\pi:X^n\to X^2$ and each $d\in D\setminus E$.
It follows that each $\pi(Z_d)$ lives in a zero-dimensional component of $D(X^2)$.
Note that if $\pi(Z_d)$ and $\pi(Z_{d'})$ live in the same zero-dimensional component of $D(X^2)$, then $\pi(Z_d)=\pi(Z_{d'})$.
Since there are only countably many irreducible components of $D(X^2)$, and only finitely many projections $\pi$, but continuum-many $d\in D\setminus E$ (as $\dim D>0$), there must exist infinitely many distinct $d_1,d_2,\dots \in D\setminus E$ with $\pi(Z_{d_i})=\pi(Z_{d_1})$ for all $i>1$ and all co-ordinate projections $\pi:X^n\to X^2$.
Applying Corollary~\ref{geom-trivial=binary} to $Z_{d_1}\subseteq X^n$,
there exists a fixed finite set of subvarieties $B_1,\dots,B_N\subseteq X^n$, such that each $Z_{d_i}$ is equal to one of the $B_j$'s.
But this contradicts the fact that the $d_i$'s are distinct and hence the $Z_{d_i}$'s are distinct (this is the uniqueness of the Douady map in the universal property for Douady spaces).
Hence, it must be the case that either $Y$ is degenerate or it lives in a zero-dimensional component of the Douady space.

By Lemma~\ref{isolated-degenerate-es}, $X$ must be essentially saturated.
\end{proof}

\bigskip

\section{Essential saturation of Inoue surfaces of type $S_M$}
\label{inouesection}

\noindent

Let us briefly recall some deformation theory of compact complex manifolds.
We suggest~\cite{kodaira} for further details and as a general reference.
A compact complex manifold $M$ is {\em rigid} if $H^1(X,T_X)=0$, where $T_X$ denotes the tangent sheaf of $X$.
Every deformation of a rigid compact complex manifold is locally trivial. More precisely: {\em Suppose $M$ is a rigid compact complex manifold.
If $\mathcal{M}\to B$ is a proper and flat surjective holomorphic map of complex varieties, and $c\in B$ is such that $\mathcal{M}_c=M$, then there exists an open neighbourhood $U$ of $c$ in $B$ such that $\mathcal{M}_U\to U$ is biholomorphic to $U\times M$ over $U$.}
Indeed, this is the classical Kodaira-Spencer deformation theory (cf. Theorem~4.6 of~\cite{kodaira})
once we observe that by flatness and the fact that $M$ is a complex manifold, the restriction of $\mathcal{M}\to B$ to some neighbourhood of $c\in B$ is a proper submersion of complex manifolds.

We will also be interested in embedded deformations, which in fact are already implicit in our discussion of Douady spaces.
Suppose $M$ is a complex submanifold of a compact complex manifold $N$.
We say that {\em $M$ is rigid in $N$} if $H^0(M,\mathcal{N}_{M/N})=0$, where $\mathcal{N}_{M/N}$ is the normal sheaf of $M$ in $N$.
In that case every deformation of $M$ in $N$ is trivial.
More precisely:
{\em Suppose $M$ is rigid in $N$.
If $B$ is a complex variety and $G$ is a complex-analytic subset of $B\times N$ such that $G\to B$ is flat and surjective, and $G_c=M$ for some $c\in B$, then $G=B\times M$.}
This follows, for example, from the fact that $H^0(M,\mathcal{N}_{M/N})$ is the tangent space to the Douady space of $N$ at the point $[M]$ corresponding to $M\subseteq N$ (cf. Proposition~1.7 of~\cite{campanapeternell94}).
Hence $[M]$ is isolated in $D(N)$ and so the Douady map $g:B\to D(N)$ must be onto the single point $[M]$.

Embedded deformations give rise to the notion of deformations of holomorphic maps that leave the domain and target fixed.
Suppose $f:M\to N$ is a holomorphic map between compact complex manifolds.
We say that {\em $f$ is rigid with respect to $M$ and $N$} if $H^0(M,f^*T_N)=0$.
It is not hard to see that this is equivalent to asking that the graph $\Gamma(f)$ is rigid in $M\times N$.
In terms of deformations of $f$ this can be formulated as follows:
{\em Suppose $f:M\to N$ is rigid with respect to $M$ and $N$.
If $B$ is a complex variety and $\Phi:B\times M\to B\times N$ is a holomorphic map over $B$ such that $\Phi_c=f$ for some $c\in B$, then $\Phi=\id_B\times f$.}

We put these facts together in the following Lemma for future use:

\begin{lemma}
\label{rigidrigid}
Let $f:M\to N$ be a holomorphic map between compact complex manifolds.
Suppose $M$ is rigid and $f$ is rigid with respect to $M$ and $N$.
Then for any flat and proper surjection of complex varieties $\mathcal{M}\to B$ with $\mathcal{M}_c=M$ for some $c\in B$, and any holomorphic map  $\Phi:\mathcal{M}\to B\times N$ over $B$ with $\Phi_c=f$, there must exist a neighbourhood $U$ of $c$ such that $\Phi_U(\mathcal{M}_U)=U\times f(M)$.
\end{lemma}

\begin{proof}
By rigidity of $M$ there exists a neighbourhood $U$ of $c$ and a biholomorphism $\sigma:U\times M\to \mathcal{M}_U$ over $U$.
We may assume that $\sigma_c=\id_M$.
So $\Phi_U\circ\sigma$ is a holomorphic map from $U\times M$ to $U\times N$ over $U$ with $(\Phi_U\circ\sigma)_c=f$.
By rigidity of $f$ with respect to $M$ and $N$,  $\Phi_U\circ\sigma=\id_U\times f$.
Hence
$$\Phi_U(\mathcal{M}_U)=\Phi_U\circ\sigma(U\times M)=(\id_U\times f)(U\times M)=U\times f(M)$$
as desired.
\end{proof}

We now specialise to the case of Inoue surfaces of type $S_M$.

\begin{lemma}
\label{rigidity}
Let $X$ be an Inoue surface of type $S_M$ and set $p_i:X\times X\to X$ to be the $i$th co-ordinate projection, for $i=1,2$.
Suppose $Y$ is an irreducible normal compact complex surface, and $f:Y\to X\times X$ is a holomorphic map with the property that $f_i:=p_i\circ f:Y\to X$ is surjective for $i=1,2$.
Then $Y$ is itself Inoue of type $S_M$.
Moreover, $f$ is rigid with respect to $Y$ and $X\times X$.
\end{lemma}

\begin{proof}
Consider the finite surjection $f_1:Y\to X$.
Let $B$ be the set of points in $Y$ at which $f_1$ is not locally biholomorphic.
Then $B$ is a complex-analytic subset of $Y$ (see~2.19 of~\cite{fischer76}).
Since $Y$ is normal and $X$ is smooth, if $B$ is nonempty then it must have dimension $1$ (see~4.2 of~\cite{fischer76}).
But as $f_1$ is finite-to-one, and $X$ contains no curves, the latter is impossible.
Hence $B$ is empty, $f_1$ is an unramified covering, and $Y$ is again Inoue of type $S_M$.
In particular, $H^0(Y,T_Y)=0$.

By the above argument $f_2:Y\to X$ is also an unramified covering.
Hence $f_i^*T_X=T_Y$ for $i=1,2$.
Since $T_{X\times X}=p_1^*T_X\oplus p_2^*T_X$, we get
$$f^*T_{X\times X}=f^*p_1^*T_X\oplus f^*p_1^*T_X=f_1^*T_X\oplus f_2^*T_X=T_Y\oplus T_Y$$
and so $H^0(Y,f^*T_{X\times X})=H^0(Y,T_Y)\oplus H^0(Y,T_Y)=0$.
So $f$ is rigid with respect to $Y$ and $X\times X$.
\end{proof}

\begin{proposition}
\label{subvarietiesofx2}
Suppose $Y\subseteq X\times X$ is an irreducible complex-analytic subset.
Then one of the following holds:
\begin{itemize}
\item[(a)]
$\dim Y=0$, or
\item[(b)]
$Y=\{a\}\times X$ or $Y=X\times\{a\}$ for some $a\in X$, or
\item[(c)]
$Y$ lives in a zero-dimensional component of the Douady space of $X\times X$.
\end{itemize}
\end{proposition}

\begin{proof}
Suppose $Y\subseteq X\times X$ is an irreducible complex-analytic subset, and let $p_i:X\times X\to X$ be the $i$th co-ordinate projection for $i=1,2$.
Since each $p_i(Y)$ is irreducible, strong minimality of $X$ implies that $p_i(Y)$ is either a point or all of $X$.
If $p_i(Y)$ is a point, then by strong minimality of $X$ one of ~(a) or~(b) must hold.
Hence we may assume $p_i(Y)=X$ for $i=1,2$.
We will show that~(c) holds.

It is clear that (c) holds if $Y = X \times X$. Let us now suppose that $Y$ is proper.
Since $p_i(Y)=X$ for $i=1,2$, $Y$ is two-dimensional by strong minimality of $X$. Let $D$ be the irreducible component of the Douady space of $X\times X$ in which $Y$ lives, and let $Z\subseteq D\times(X\times X)$ be the restriction of the universal family to $D$.
Let $\widetilde{Z}\to Z$ be the normalisation of $Z$.
We have the following situation
$$\xymatrix{
\widetilde{Z}\ar[r]^{\pi}\ar[dr] & Z\ar[r]^{\subseteq \ \ \ \ \ \ \ \ }\ar[d] &D\times (X\times X)\ar[dl]\\
& D
}$$
There exists a proper complex-analytic subset $E\subset D$ such that for all $d\in D\setminus E$,
\begin{itemize}
\item[(i)]
$Z_d\subset X\times X$ is a reduced and irreducible surface,
\item[(ii)]
$p_i(Z_d)=X$ for $i=1,2$,
\item[(iii)]
$\pi_d:\widetilde{Z}_d\to Z_d$ is the normalisation of $Z_d$, and
\item[(iv)]
$\widetilde{Z}\to D$ is flat outside of $E$.
\end{itemize}
Indeed, we can choose $E$ to satisfy~(i) though~(iv) because $Z\to D$ is flat and proper.
For~(i), use the fact that $Z\to D$ has one reduced and irreducible two-dimensional fibre (namely $Y$) and hence, by flatness, its general fibres are so.
For~(ii), use the fact that none of the $Z_d$'s are degenerate since $Y$ was not (cf. Lemma~\ref{degeneratecomponents}), and we have already seen that this implies its projections to $X$ are surjections.
To find $E$ satisfying~(iii) note that the general fibre of $\widetilde{Z}\to D$ is normal (Th\'eor\`eme~2 of~\cite{banica79}) and $\pi$ restricted to the general fibre is again a finite map that is biholomorphic outside a proper complex-analytic set.
Finally, we can find $E$ satisfying~(iv) since by~\cite{frisch67} every proper holomorphic map is flat outside a proper complex-analytic set.

Now fix $d_0\in D\setminus E$.
Then $\widetilde{Z}_{d_0}$ is an irreducible normal compact complex surface such that $\pi_{d_0}:\widetilde{Z}_{d_0}\to X\times X$ composed with the projections to $X$ are surjective.
By Lemma~\ref{rigidity}, $\widetilde{Z}_{d_0}$ is itself an Inoue surface of type $S_M$ -- and hence rigid -- and the map $\pi_{d_0}$ is rigid with respect to $\widetilde{Z}_{d_0}$ and $X\times X$.
By Lemma~\ref{rigidrigid}, there exists an open neighbourhood $U$ of $d_0$ in $D\setminus E$, such that $\pi_U(\widetilde{Z}_U)=U\times Z_{d_0}$.
Hence $Z_U=U\times Z_{d_0}$.
The universal property of the Douady space implies that $U=\{d_0\}$.
But as $U$ was open in the irreducible $D$, this means that $D=\{d_0\}$.
We have shown that $Y$ lives in a zero-dimensional component of the Douady space of $X\times X$, as desired.
\end{proof}

\begin{corollary}
Inoue surfaces of type $S_M$ are essentially saturated but not of K\"ahler-type.
\end{corollary}

\begin{proof}
By Fact~\ref{inouefacts}, Inoue surfaces of type $S_M$ are strongly minimal compact complex varieties not of K\"ahler-type.
By Fact~\ref{trichotomy} they are trivial.
Proposition~\ref{subvarietiesofx2} tells us that the irreducible subvarieties of $2$-space are either degenerate or live in a  zero-dimensional component of the Douady space.
By Proposition~\ref{reductionto2} this is then true of $n$-space for all $n\geq 0$, and these surfaces are essentially saturated.
\end{proof}

\begin{remark}
One might be tempted to try and generalise the above corollary by observing that our argument goes through for any compact complex surface satisfying the four properties of Fact~\ref{inouefacts}.
However, these properties actually characterise Inoue surfaces of type $S_M$.
Indeed, from Kodaira's classification of compact complex surfaces we see that a non-K\"ahler surface $X$ without curves must have $b_1(X)=1$.
It will also have $b_2(X)=0$ by Riemann-Roch if one imposes
the condition $H^1(X,T_X)=0$.
But, as mentioned before, such surfaces are completely classified (cf.~\cite{inoue74, LYZ94, Tel94}) and they belong to one of Inoue's classes.
Since the surfaces of type $S^{(+)}$ have $H^1(X,T_X)=1$, we are left with only the surfaces of type $S_M$ and $S^{(-)}$.
Among them only those of type $S_M$ satisfy the fourth property, since any surface of type  $S^{(-)}$ admits a double cover of type $S^{(+)}$, see \cite{inoue74}.
\end{remark}

We conclude by pointing out that among the known strongly minimal compact complex surfaces, the surfaces of type $S_M$ are the only non-K\"ahler essentially saturated examples.
Indeed, since the only known non-K\"ahler compact complex surfaces without curves are the Inoue surfaces, this follows from:

\begin{proposition}
The other Inoue surfaces, those of type $S^{(+)}$ and $S^{(-)}$, are not essentially saturated.
\end{proposition}

\begin{proof}
Note that by Fact~\ref{trichotomy}, these strongly minimal surfaces are also trivial.

Let $X$ be an Inoue surface of type $S^{(+)}$.
The universal cover of $S^{(+)}$ (and indeed of all the Inoue surfaces) is $\mathbb{H}\times\mathbb{C}$, the product of the upper half plane with the complex plane.
From Inoue's construction of $S^{(+)}$ it is evident that translation on the second co-ordinate induced a non-trivial action of $(\mathbb{C},+)$ on $X$ (see equation~(18) of~\cite{inoue74}).
We thus obtain an infinite analytic family of automorphisms of $X$ paramaterised by $\mathbb{C}$, which by considering graphs, can be viewed as living in the irreducible component $D$ of $D(X\times X)$ that contains the diagonal $A\subset X\times X$.
This already implies that $D$ is not compact, since no trivial strongly minimal compact complex manifold can have an infinite definable family of automorphisms.
However, the non-compactness of $D$ can be seen more directly, without using triviality of $X$, as follows:
Assume $D$ is compact and let $Z\subset D\times X\times X$ be the universal family over $D$.
Since $\dim  H^0(X,T_X)=1$, and $\mathcal{N}_{A/X\times X}\cong T_X$, we know that $\dim D=1$.
Fixing $a\in X$, the subspace $Z\cap (D\times\{a\}\times X)$
is thus one-dimensional and contains $\{(g,a,ga): g\in \mathbb{C}\}$.
Thus its projection on $X$ is a complex-analytic subset containing the orbit of $a$ under the action of $(\mathbb C,+)$, and therefore cannot be zero dimensional.
This contradicts the fact that $X$ has no curves.

Let now $Y$ be an Inoue surface of type $S^{(-)}$.
As mentioned before, there exists an Inoue surface $X$ of type $S^{(+)}$ which is a double cover of $Y$.
Now one can finish by showing that a finite cover of an essentially saturated compact complex variety is again essentially saturated.
But in this case the argument is easier:
The action of $(\mathbb{C},+)$ on $X$ described above will induce an infinite analytic family of subvarieties of $Y\times Y$ which project in a finite-to-one manner onto each component.
Arguing as above one shows that this family of subspaces of $Y\times Y$ cannot be compactified in $D( Y\times Y)$.
(Alternatively, essential saturation of $Y$ would imply the existence of an infinite definable family of finite-to-finite correspondences, which is also ruled out by the triviality of $Y$.)
\end{proof}


\end{document}